\documentclass{amsart}
\usepackage{amssymb} 
\usepackage{amsthm, multirow, subfigure}
\usepackage[all]{xy}
\usepackage{graphicx}

\newtheorem{theorem}{Theorem}
\newtheorem*{corollary*}{Corollary}
\newtheorem{lemma}{Lemma}[section]
\newtheorem{proposition}{Proposition}[section]
\newtheorem{corollary}{Corollary}[section]
\theoremstyle{definition}

\theoremstyle{remark}

\title{Substitutions and $\frac{1}{2}$-Discrepancy of $\{n \theta + x\}$ \rm{II}}
\author{David Ralston}
\email{ralston.david.s@gmail.com}
\address{Ben Gurion University, Department of Mathematics \\ POB 653 \\ Beer Sheva 84105 \\ ISRAEL}
\date{\today}
\begin{document}
\maketitle
\section{Introduction}
Given an irrational $\theta$ and some $x \in [0,1)$, the \textit{$1/2$-discrepancy sums of $x$} are given by
\[S_i(x) = \sum_{j=0}^{i-1} \left(\chi_{[0,1/2)}-\chi_{[1/2,1)}\right)(x+j \theta),\] where all addition is taken modulo one.  We are interested in studying the asymptotic growth of $S_i(x)$ for generic $\theta$.  However, as $S_i(x)$ are not monotone, it is convenient to define
\begin{equation}\label{eqn - rho}
\rho_n(x) = 1 + \max\left\{S_i(x): i=1,\ldots,n\right\} - \min\left\{S_i(x): i=1,\ldots,n\right\}.
\end{equation}

In \cite{ralston}, it was shown that there is a natural renormalization procedure for studying this sequence; we will briefly reintroduce the necessary notation in \S\ref{section - refresher}, but familiarity with the existing techniques will be helpful.  At the center of this renormalization procedure is a map $g:S^1 \rightarrow S^1$, closely related to the Guass map.  While that work was concerned with developing the renormalization procedure in order to construct specific behaviors, we will begin our study of \textit{generic} growth rates by showing first:

\begin{theorem}\label{theorem - g ergodic}
There is a unique measure $\mu_g$ on the circle which is mutually absolutely continuous with respect to Lebesgue measure and is preserved by $g$.  Both Radon-Nikodym derivatives $dx/d\mu_g$ and $d\mu_g/dx$ are essentially bounded.  Furthermore, the system $\{S^1,\mu_g, g\}$ is exponentially CF-mixing.
\end{theorem}

If $\rho_n(x) \in o(b_n)$ for \textit{every} $x \in [0,1)$, then we will simply wrote $\rho_n \in o(b_n)$.  Similarly, we will write $\rho_n \notin o(b_n)$ if for \textit{every} $x$ we have $\rho_n(x) \notin o(b_n)$.  Note that trivially $S_n(x) \in o(b_n)$ for every $x$ if and only if $\rho_n \in o(b_n)$, and similarly for $S_n(x) \notin o(b_n)$.

A function $f: \mathbb{R}^+ \rightarrow \mathbb{R}$ is said to be \textit{regularly varying} if for any $C>0$ we have $f(Cx) \sim f(x)$.

\begin{theorem}\label{theorem - main result}
Suppose that $f(x)$ is a continuous nondecreasing function which is regularly varying, defined for all $x \geq C$.  Define the sequence $\{F_n\}$ by
\begin{equation}\label{eqn - defining B} F(t) = \int_C^t f(x)dx.\end{equation}
Then for almost every $\theta$, either $\rho_n \in o(F(\log n))$ or $\rho_n \notin o(f(\log n))$ according to whether
\[\int_C^{\infty}f(x)dx < \infty \quad \textrm{or} \quad \int_C^{\infty}f(x)dx=\infty.\]
\end{theorem}

\section{Existing Notation and Prior Results}\label{section - refresher}

All notation is consistent with \cite{ralston}; whenever the same objects are defined here as in that work, the same notation will be used.

We use standard continued fraction notation, and as $\theta \in (0,1)$ without loss of generality, we omit the integer part and write for $a_i \in \mathbb{N}$
\[\theta=\frac{1}{a_1+\cfrac{1}{a_2+\ddots}}=[a_1,a_2,\ldots].\]  The partial quotients may be written as $a_i(\theta)$ when $\theta$ is not immediately clear from context.  The Gauss map acts as the one-sided shift on the sequence of partial quotients:
\[\gamma(\theta) = \frac{1}{\theta}-a_1, \quad \gamma([a_1,a_2,\ldots])=[a_2,a_3,\ldots].\]
The Gauss map preserves a unique probability measure which is mutually absolutely continuous with respect to Lebesgue measure, and it is exponentially CF-mixing with respect to this measure.  We denote this measure by $\mu_{\gamma}$, and it is given by
\[\mu_{\gamma}(A) = \frac{1}{\log 2} \int_A \frac{dx}{1+x}.\]

The function $g$ is given by
\[g([a_1,a_2,a_3,\ldots])= \begin{cases}[a_2+1,a_3,\ldots] = 1 - \theta & (a_1=1)\\ [1,a_2,a_3,\ldots] = \frac{1}{1+\gamma(\theta) } & (a_1 = 1 \bmod 2, \, \neq 1)\\ [a_3,a_4,\ldots]=\gamma^2(\theta) & (a_1=0\bmod 2).\end{cases}\]

We define
\[E(x) = \max\{n \leq x: n \in \mathbb{Z}, \, n=0 \bmod 2\},\] and we then let
\[\theta_n = g^n(\theta), \quad \delta_n = 1-E(a_1(\theta_n))\theta_n.\]

If we assume that $\theta=\theta_0 <1/2$, then we may label the intervals
\[A=[0,1/2), \quad B=[1/2,1-\theta), \quad C=[1-\theta,1).\]

We omit a discussion of standard terminology in substitutions systems (see \cite[\S 2]{ralston}), except to point out that most techniques are carried out in the compact but disconnected symbol space in which $A$, $B$ and $C$ are all compact.  \textit{This distinction can only cause problems in regards to the endpoints of the intervals, and any orbit can include at most two of the endpoints; we do not in general distinguish between $S^1$ and the associated symbol space.}  Define the substitutions $\sigma_n=\sigma(\theta_n)$ according to Table \ref{table - substitutions}.  For convenience denote
\[\sigma^{(n)}=\sigma_0 \circ \sigma_1 \circ \cdots \circ \sigma_{n-1}.\]

\begin{table}[hbt]
\begin{tabular}{| l | l |}
\hline Case & Substitution\\
\hline \hline \multirow{3}{*}{$a_1=2k, \, a_3 \neq 1$} &$A \rightarrow (A^{k+1}B^{k-1}C) (A^{k}B^{k-1}C)^{a_2-1}$\\
& $ B \rightarrow (A^{k}B^{k}C) (A^{k}B^{k-1}C)^{a_2-1}$\\
& $C \rightarrow (A^{k}B^{k}C) (A^{k}B^{k-1}C)^{a_2}$\\
\hline \multirow{3}{*}{$a_1=2k, \, a_3=1$} & $A \rightarrow (A^{k}B^{k}C) (A^{k}B^{k-1}C)^{a_2}$\\
& $B \rightarrow (A^{k+1}B^{k-1}C)(A^{k}B^{k-1}C)^{a_2}$\\
& $C \rightarrow (A^{k+1}B^{k-1}C)(A^{k}B^{k-1}C)^{a_2-1}$\\
\hline \multirow{3}{*}{$a_1=2k+1$}& $A \rightarrow A^{k}B^{k}C$\\
&$B \rightarrow A^{k+1}B^{k-1}C$\\
&$C \rightarrow A$\\
\hline \multirow{3}{*}{$a_1=1$} & $A \rightarrow A$\\
&$ B \rightarrow B$\\
&$C \rightarrow C$\\
\hline
\end{tabular}
\caption{The substitution $\sigma$ as a function of $\theta$.}
\label{table - substitutions}
\end{table}

Under our labeling of $A$, $B$ and $C$, we may without confusion define for finite words $\Omega$ in the alphabet $\{A,B,C\}$
\[S(\Omega) = \#A - \#B - \#C,\] and denoting by $\Omega(n)$ the first $n$ letters in $\Omega$ (for $0<n\leq |\Omega|$, the length of $\Omega$), we may similarly define
\[\rho(\Omega) = 1 + \max\{S(\Omega(n)): n=1,2,\ldots,|\Omega|\} - \min\{S(\Omega(n)): n=1,2,\ldots,|\Omega|\}.\]

By \cite[Theorem 1.1]{ralston}, there is a special point $x(\theta)$ for which the sequence of words
\begin{equation}\label{eqn - omega n}\Omega_n=\sigma^{(n)}(A)\end{equation} \textit{correctly encodes the orbit} of $x(\theta)$ under rotation by $\theta$ with at most two errors.  That is, with two possible exceptions, the sequence of letters in $\Omega_n$ correctly specifies in which interval $x+i\theta$ lies.  Combined with \cite[Proposition 5.2]{ralston}, we have
\begin{equation}\label{eqn - computing spread of special x}\rho(\Omega_n) = \frac{1}{2}\sum_{i=0}^n E(a_1(\theta_n)) + \xi,\end{equation} where $|\xi| \leq 5$ for all $\theta$.  As $\{S_n(x)\}$ must be unbounded for every $x$ \cite{MR0209253}, \textit{we will ignore the bounded term $\xi$ from this point forward}.

The words $\Omega_n$ represent the encoding of the first return of $x(\theta)$ to an interval $\tilde{I}_n$.  The intervals $\tilde{I}_n$ are nested and compact (in the associated symbol space), and of length
\[|\tilde{I}_n|=\delta_0 \cdot \delta_1 \cdots \delta_{n-1}.\]
The point $x(\theta)$ is the intersection of all $\tilde{I}_n$; we point out that for all irrational $\theta$
\begin{equation}\label{eqn - exponential decay of I_n}
\liminf_{n \rightarrow \infty} \left|\frac{\log (\delta_0 \cdots \delta_{n-1})}{n} \right| \geq \frac{\log 2}{2} > 0,
\end{equation}
which is to say that the length of $\tilde{I}_n$ always decays at least exponentially fast.  This observation is direct from the fact that no two consecutive $\theta_i$ may both be larger than $1/2$, and for $\theta_i<1/2$ we have $\delta_i < 1/2$.

While not explicitly stated in the proof of \cite[Theorem 1.1]{ralston}, the following is immediate from the techniques of the proof, especially \cite[Proposition 4.1]{ralston}:
\begin{corollary*}
For each $y \in [0,1)$ and every $n$, there are two words $\Upsilon_0(y,n)$, $\Upsilon_1(y,n)$, where
\[\Upsilon_1(y,n) \in \{\sigma^{(n)}(A),\sigma^{(n)}(B),\sigma^{(n)}(C)\},\] and $\Upsilon_0(y,n)$ is either a proper right factor of one of these words or empty, and the concatenated word
\begin{equation}\label{eqn - words for orbit of any y}\Omega_n(y) = \Upsilon_0(y,n) \Upsilon_1(y,n) \end{equation} correctly encodes the orbit of $y$, except for at most two errors.
\begin{proof}
By \cite[Theorem 1.1]{ralston}, the induced transformation on each $\tilde{I}_n$ is rotation by $\theta_n$, and we may therefore label intervals $A$, $B$ and $C$ within $\tilde{I}_n$ according to this new rotation $\theta_n$; the specific labeling is outlined in that proof.  Let $\Upsilon_0(y,n)$ be the word which encodes the orbit of $y$ (in the disconnected symbol space corresponding to the original transformation) through its return to $\tilde{I}_n$ (if $y \in \tilde{I}_n$, then let $\Upsilon_0(y,n)$ be empty), and let $\Upsilon_1(y,n)$ be $\sigma^{(n)}(*)$, where $*$ represents which of the intervals $A$, $B$, $C$ contains this first point in the orbit of $y$ (under the correct labeling of intervals in the induced system on $\tilde{I}_n$).  As the collection \[\{\sigma^{(n)}(A), \sigma^{(n)}(B), \sigma^{(n)}(C)\}\] encodes all possible orbits of points in $\tilde{I}_n$ until their return to $\tilde{I}_n$ and the original rotation by $\theta$ is minimal, $\Upsilon_0$ must be a right factor of one of these words.
\end{proof}
\end{corollary*}

Define the matrices $M_n=M(\theta_n)$ according to Table \ref{table - matrices for return times} (the eigenvalues are readily computable).  For convenience, denote
\[M^{(n)}=M_{n-1}\cdot M_{n-2} \cdots M_0,\] so by \cite[Lemma 5.4]{ralston}, we have both that $|\sigma^{(n)}(A)|=|\sigma^{(n)}(B)|$ and
\begin{equation}\label{eqn - lengths of words} M^{(n)} \left[ \begin{array}{c} 1 \\ 1\end{array}\right] = \left[ \begin{array}{c} |\sigma^{(n)}(A)| \\ |\sigma^{(n)}(C)| \end{array} \right].\end{equation}

\begin{table}[hbt]
\begin{tabular}{| c | c | c|}
\hline Case & $M(\theta)$ & Eigenvalues\\
\hline \hline $a_1=2k$, $a_3 \neq 1$ & $\left[ \begin{array}{c c}(a_1-1)a_2+1 &a_2 \\ (a_1-1)a_2+a_1&a_2+1 \end{array}\right]$ & $ka_2+1 \pm \sqrt{ka_2(ka_2+2)}$\\
\hline $a_1=2k$, $a_3 = 1$ & $\left[ \begin{array}{c c} (a_1-1)a_2+a_1 &a_2+1 \\(a_1-1)a_2+1 &a_2 \end{array}\right]$ & $ka_2+1 \pm \sqrt{ka_2(ka_2+2)}$\\
\hline $a_1=2k+1$ & $\left[ \begin{array}{c c} a_1-1& 1 \\ 1& 0 \end{array}\right]$ & $k \pm \sqrt{k^2+1}$\\
\hline $a_1=1$ & $\left[ \begin{array}{c c} 1& 0 \\ 0 &1\end{array}\right]$ & $1$\\
\hline
\end{tabular}
\caption{The matrices $M(\theta)$ used to determine the lengths of $\sigma^{(n)}(A)$, $\sigma^{(n)}(C)$.  Note that $2k=E(a_1)$.}\label{table - matrices for return times}
\end{table}

\section{Proof of Theorem \ref{theorem - g ergodic}}
Note that $g^{-1}(1/2,1) \subset (0,1/2)$, so any probability measure $\mu$ which is preserved by $g$ must have
$\mu(1/2,1) \leq 1/2$: the Gauss measure $\mu_{\gamma}$ is not preserved by $g$.

Define the following collection of open intervals (for $n,m,k=1,2,\ldots$):
\[ \left( \frac{1}{2},1\right), \quad \left(\frac{1}{2k+2},\frac{1}{2k+1}\right), \quad \left(\frac{m}{2nm+1},\frac{m+1}{2n(m+1)+1}\right).\]
The middle intervals are those $\theta$ whose continued fraction expansion begins with an odd number (except one); we will refer to the collection of all such intervals as \textit{odd intervals}.  The right-most intervals are those $\theta$ whose continued fraction expansion begins with the pair $[2n,m,\ldots]$, and we will refer to the collection of all such intervals as \textit{even intervals}.  The collection of odd and even intervals, together with $(1/2,1)$ we denote $\mathcal{C}$, and note that $\mathcal{C}$ covers $S^1$ except for a countable set of points.

\begin{lemma}\label{lemma - invariant density properties}
We have each of the following:
\renewcommand{\labelenumi}{\Roman{enumi}.}
\renewcommand{\theenumi}{(\Roman{enumi})}
\begin{enumerate}
\item $\mathcal{C}$ is a Markov partition for $g$. \label{item - markov}\\
\item For $c_1,c_2 \in \mathcal{C}$, if $g(c_1)\cap c_2 \neq \emptyset$, then $c_2 \subset g(c_1)$. \label{item - covering}\\
\item There is a $k$ such that $(0,1)\subset g^k(c)$ for each $c \in \mathcal{C}$. \label{item - covering 2}\\
\item The map $g$ is monotone and $1:1$ on each $c \in \mathcal{C}$, and extends to the closure of each $c$ to a $C^2$ function. \label{item - C2}\\
\item $g$ is expansive; $\exists k$ and $d>1$ such that $|(g^k)'(x)|\geq d$ for almost all $x$. \label{item - expansive}\\
\item $g$ has the \textit{Renyi} (or \textit{strong distortion}) property: $\exists d$ such that for every $c \in \mathcal{C}$ \[\sup_{x \in c} \left(\frac{|g''(x)|}{(g'(x))^2}\right)<d.\]\label{item - renyi}
\item The endpoints of all $c \in \mathcal{C}$ map (via the $C^2$ extension from item \ref{item - C2}) to a finite set.\label{item - endpoints}
\end{enumerate}
\begin{proof}
The restriction of $g$ to each $c \in \mathcal{C}$ is invertible, and we have chosen $\mathcal{C}$ to generate the Borel $\sigma$-algebra under $g^{-1}$; if $\theta_1 \neq \theta_2$, then there is some minimal index $i$ such that $a_i(\theta_1)\neq a_i(\theta_2)$, from which one sees that there is some $k$ such that either (see for example \cite[Eqn. 20]{ralston})
\[a_1(g^k \theta_1) \neq a_1(g^k \theta_2), \quad \textrm{or}\, a_1(g^k \theta_1)=a_1(g^k \theta_2)=0\bmod 2, \quad a_2(g^k \theta_2) \neq a_2(g^k \theta_2).\]  So item \ref{item - markov} is shown.

The related items \ref{item - covering} and \ref{item - covering 2} can both be shown directly using the fact that $g$ maps odd intervals to $(1/2,1)$, $g$ maps $(1/2,1)$ to $(0,1/2)$, and $g$ maps even intervals to $(0,1)$.

On the even intervals, $g=\gamma^2$ is the square of the Gauss map, and the even intervals are members of the standard Markov partition for $\gamma^2$, which is well-known to have all of these properties.  Furthermore, on $(1/2,1)$ we have $g(\theta)=1-\theta<1/2$.  Using the chain rule, then, we need only establish the remaining items for the odd intervals.  So, let us consider the odd interval for a fixed $k$.  The reader may verify that for $\theta=[2k+1,\ldots]$ we have
\[\frac{1}{2k+1}<1-2k\theta<\frac{1}{k+1}.\]
Using this inequality, one may show:
\begin{align*}
g(\theta) &= \frac{1}{1+\gamma(\theta)}=\frac{\theta}{1-2k \theta},\\
g'(\theta) &= \frac{1}{(1-2k\theta)^2}\geq (k+1)^2,\\
g''(\theta) &=\frac{4k}{(1-2k\theta)^3}\leq 4k(2k+1)^3.\\
\end{align*}
Item \ref{item - C2} is immediate from the first line (for $\theta$ in this odd interval, $\theta<(2k+1)^{-1}$, so there is no asymptote), and item \ref{item - expansive} from the second (recall that $k \geq 1$).  Item \ref{item - renyi} follows from the bounds on $g'$ and $g''$.
\end{proof}
\end{lemma}

The literature surrounding the existence of invariant measures for Markov maps of the interval is vast.  We have arranged Lemma \ref{lemma - invariant density properties} to match the statement which appears (as an unnumbered theorem) in \cite{MR556585}.  This theorem gives the existence of a unique probability density $\mu_g$, supported on all of $[0,1)$, which is invariant under $g$ and continuous with respect to Lebesgue measure, with essentially bounded Radon-Nikodym derivative.  Therefore $\{S^1,\mu_g, g\}$ is ergodic.  Following the trail of references in that work back to \cite{MR0486431} actually says more, however, that the system $\{S^1,\mu_g, g\}$ is weakly Bernoulli (and therefore exact).  Item \ref{item - covering 2} gives that $g$ is topologically mixing (any open set contains an open subinterval specified by an initial finite string of partial quotients), so by \cite[Corollary 4.7.8]{MR1450400}, our system is exponentially CF-mixing.  That the Radon-Nikodym derivative $d\mu_g/dx$ is bounded away from zero (not just bounded) is given as a remark in the third paragraph of \cite[\S 4.7]{MR1450400}.

The following Khinchin-like characterization follows immediately via a standard shrinking-target result (or the Borel-Cantelli Lemma) as $g$ is mixing and $\mu_g \sim \mu_{\gamma}$, the measure preserved by the Gauss map (as both are mutually absolutely continuous with respect to Lebesgue measure):
\begin{corollary}\label{corollary - khinchin like}
Let $\{b_i\}$ be a sequence of positive real numbers for $i=0,1,\ldots$.  Then the inequality
\[a_1(g^n \theta) > b_n\] is satisfied almost surely infinitely many times or only finitely many times according to whether the series
\[\sum_{i=0}^{\infty} \frac{1}{b_i}\] diverges or not.
\end{corollary}

\section{Proof of Theorem \ref{theorem - main result}}

For a square real-valued matrix $M$, let $\|M\|$ be the largest eigenvalue, and for a real-valued column vector $u$, let $\|u\|$ be the largest element (in absolute value).
\begin{lemma}
With $M(\theta)$ as given by Table \ref{table - matrices for return times}, we have both $\log \|M(\theta)\|$ and $\log \|M^{-1}(\theta)\|$ in $L^1(X, \mu_g)$.
\begin{proof}
As the two-by-two matrices $M(\theta)$ all have $| \det{M(\theta)} |=1$, we have $\|M(\theta)\|=\|M^{-1}(\theta)\|$.  As $\mu_{g}$ is mutually absolutely continuous with Lebesgue measure, it therefore suffices to show that
\[\int_{0}^1 \log \|M(\theta)\| d\theta < \infty.\]

For $\theta \in (1/2,1)$, we have $\|M(\theta)\|=1$.  If $a_1(\theta)=2k+1$ for $k \neq 1$, then
\[\log \|M(\theta)\| = \log | k + \sqrt{k^2+1}| < \log(2k+1).\]
On the other hand, for $a_1(\theta)=2n$ and $a_2(\theta)=m$ we have
\[\log \|M(\theta)\| = \log | nm+1 + \sqrt{nm(nm+2)}| \leq \log(2nm+2) .\]
So we may therefore compute:
\begin{align*}
\int_0^1 \log \|M(\theta)\|d \theta &\leq \sum_{k=1}^{\infty} \int_{\frac{1}{2k+2}}^{\frac{1}{2k+1}} \log (2k+1) d\theta\\
&\quad + \sum_{n,m=1}^{\infty} \int_{\frac{m}{2nm+1}}^{\frac{m+1}{2n(m+1)+1}}\log(2nm+2)d \theta\\
&= \sum_{k=1}^{\infty} \frac{\log(2k+1)}{(2k+1)(2k+2)}+ \sum_{n,m=1}^{\infty} \frac{\log(2nm+2)}{(2nm+1)(2n(m+1)+1)},
\end{align*}
and the summability of both series is direct.
\end{proof}
\end{lemma}

By the Oseledec ergodic theorem, then, for almost every $\theta$ we have
\begin{equation}\label{eqn - used oseledec}
\lim_{n \rightarrow \infty} \frac{1}{n} \log \left( M^{(n)} \left[ \begin{array}{c} 1 \\ 1 \end{array}\right]\right) = \lambda
\end{equation}
for some $\lambda<\infty$.

\begin{proposition}\label{proposition - control word growth}
For almost every $\theta$, for every $\epsilon > 0$ we have for sufficiently large $n$ (recall \eqref{eqn - omega n})
\[ (\lambda - \epsilon)^{n-3} \leq \left| \Omega_n \right| \leq (\lambda+\epsilon)^{n}.\]
\begin{proof}
The upper inequality is direct in light of the previous remarks and \eqref{eqn - lengths of words}.  For the lower inequality we must establish
\[\min \{ |\sigma^{(n)}(A)|, |\sigma^{(n)}(C)| \} \geq \max \{ | \sigma^{(n-3)}(A)|, |\sigma^{(n-3)}(C)|\}.\]
The proof may be accomplished through an exhaustive case-by-case analysis of different possible forms for the matrices $M_{n-1}$, $M_{n-2}$ and $M_{n-3}$.  The situation is easiest in the case that $a_1(\theta_{n-1})=0 \bmod 2$.  For example, suppose that
\[M_{n-1} = \left[ \begin{array}{c c } (2k-1)m+1 & m \\ (2k-1)m+2k & m+1\end{array}\right], \quad M^{(n-1)}u=\left[\begin{array}{c} A \\ B \end{array}\right].\]
Then we have
\[M^{(n)}u = \left[ \begin{array}{c} C \\ D \end{array}\right] = \left[ \begin{array}{c} \left((2k-1)m+1\right)A+mB\\ \left((2k-1)m+2k\right)A+(m+1)B\end{array}\right],\] and we clearly have
\[\min\{C, D\} \geq \max\{A, B\}.\]
The matrices $M(\theta)$ for $a_1(\theta)=1 \bmod 2$ are less trivial, but the composition of two such matrices is seen to have the desired property. Since two such matrices can only occur separated by an identity matrix, and it is possible that $M_{n-1}$ was the identity matrix, the $n-3$ in the exponent is sufficient.
\end{proof}
\end{proposition}
\begin{lemma}
$\lambda>1$.
\begin{proof}
See \eqref{eqn - exponential decay of I_n}.  As the length of $\tilde{I}_n$ decays at least exponentially fast, the return time of any point in $\tilde{I}_n$ to itself increases exponentially fast.  As the entries of $M^{(n)}u$ are the two return times of points in $\tilde{I}_n$ to itself \eqref{eqn - lengths of words}, we must have $\lambda \geq \sqrt{2}$.
\end{proof}
\end{lemma}
Now let $f(x)$ be continuous, nondecreasing and regularly varying, and define $F(t)$ as in the introduction \eqref{eqn - defining B}.  We proceed now under the assumption that $f(x)$ is not integrable.  From Corollary \ref{corollary - khinchin like}, it follows that for generic $\theta$ and infinitely many $n$, we have $f(n) < (1/2)E(a_1(\theta_n))$.
So by \eqref{eqn - computing spread of special x},
\[ \rho(\Omega_n)  > f(n).\]
By Proposition \ref{proposition - control word growth}, then, we have infinitely many times $N$ such that
\[\rho_N(x(\theta)) \geq f(C \log N),\] where $C = \log (\lambda + \epsilon)$, so as $f$ was regularly varying, $\rho_n(x(\theta)) \notin o(f(\log n))$.

Similarly, if $f$ is integrable, we eventually have for generic $\theta$ that $(1/2)E(a_1(\theta_n))<f(n)$, so up to a bounded difference which we ignore, we may say that for generic $\theta$ we eventually have
\[\frac{1}{2}\sum_{i=0}^{n-1}E(a_1(\theta_i))<F(n),\]
from which it follows (by considering all $|\Omega_{n-1}|\leq N \leq |\Omega_{n}|$ and using regularity of $f$ and Proposition \ref{proposition - control word growth} as before) that $\rho_n(x(\theta)) \in O(F(\log n))$.  As we may multiply $f$ by any $\epsilon>0$ without affecting summability, then, we have $\rho_n(x(\theta)) \in o(F(\log n))$.

It remains to show that the behavior of \textit{any} $\rho_n(y)$ may be accurately considered through the sequence $\rho_n(x(\theta))$.
Note that regardless of the sequence of substitutions $\sigma_i$, we have (refer to Table \ref{table - substitutions} and \cite[Proposition 5.1]{ralston})
\[\rho(\sigma^{(n-1)}(A)) \leq \rho(\sigma^n(*)) \leq \rho(\sigma^n (A)),\] where $* \in \{A,B,C\}$.  Finally, it is clear that for any two words $\nu_1$ and $\nu_2$, we have
\[\rho(\nu_1 \nu_2) \leq \rho(\nu_1)+\rho(\nu_2).\]  Altogether, then, recall \eqref{eqn - words for orbit of any y}, to see that for \textit{any} $y$, for \textit{every} $\theta$ and \textit{any} $n \geq 1$:
\[\rho(\Omega_{n-1}) \leq \rho(\Omega_n(y)) \leq 2 \rho(\Omega_n).\]

This small uncertainty in the estimation gives for any $\epsilon>0$ and sufficiently large $n$
\[(\lambda-\epsilon)^{n-4} \leq |\Omega_{n-1}| \leq |\Omega_n(y)| \leq 2\max\{|\sigma^{(n)}(A)|,|\sigma^{(n)}(C)|\} \leq 2(\lambda+\epsilon)^{n}.\]
All previous arguments using regularity of $\{b_n\}$ still apply, then, extending the existing arguments to all points for generic $\theta$ and completing the proof of Theorem \ref{theorem - main result}.

\section{Concluding Remarks}\label{section - concluding remarks}

We present a pair of simple computations to show the relatively tight control we may generically impose upon $\rho_n(y)$ for any $y$ and for generic $\theta$.  A classical application of the Denjoy-Koksma inequality is to show that for $\theta$ having bounded partial quotients, for every $x$ we have $\rho_n(x) \in O(\log n)$ (see for example \cite[\S2.1]{conze}).  It was shown in \cite[Theorem 1.4]{ralston} that in fact $\rho_n(x) \sim \log(n)$ for such $\theta$.  However, by setting $f(x)=x \log x \cdots \log^{(k-1)}x$, the product of the first $k$ iterated logarithms (starting at index zero, and defined for sufficiently large $x$), we see that for generic $\theta$ we have
\[\rho_n \notin o(\log n \cdots \log^{(k)} n).\]

Another classical application of the Denjoy-Koksma inequality is that for generic $\theta$, for every $x$ we have $\rho_n(x) \in o(n^{\epsilon})$ for every $\epsilon>0$ (also in \cite[\S2.1]{conze}).  To improve upon this bound, consider
\[f(x) = x \log x \cdots \log^{(k-2)}x (\log^{(k-1)} x)^{1+\epsilon}\]
for sufficiently large $x$.  Using the straightforward fact (one may estimate the integral with the corresponding sum over integers) that
\[\int_C^t x \log x \cdots \log^{(k)}(x) \sim t \log t \cdots \log^{(k)}t\] (setting $C$ sufficiently large so that all terms are defined), we have
\begin{align*}
F(\log n) &< \left(\log^{(k-1)}(\log n)\right)^{1+\epsilon} \left( \sum_{i=1}^{\log n} i \log i \cdots \log^{(k-2)}(n) \right) \\
&\sim \log n \cdots \log^{(k-1)}(n) (\log^{(k)}(n))^{1+\epsilon},\end{align*} so generically our discrepancy sums grow slower than all such functions.

One might hope that there could be some sequence such that generically $\rho_n \sim b_n$.  Such a quest would be a fool's errand: the sequence $\rho(\Omega_n)$ is given by a partial ergodic sum of the nonintegrable nonnegative function $E(a_1 \theta)$, and the growth rate is therefore seen to be similar to the ergodic sums of simply $a_1(\theta)$.  While the gap $g$ and the Gauss map $\gamma$ are not identical and do not preserve the same measure, they are both exponentially CF-mixing with respect to mutually absolutely continuous measures.

The partial sums $a_1+a_2+\cdots+a_n$ (ergodic sums of the function $a_1(\theta)$ under the action of $\gamma$) almost surely do not have a strong law of large numbers (a result known already to Khinchin \cite{MR1556899}).  However, that example \textit{does} admit a weaker `trimmed' law of large numbers \cite{diamond-vaaler}: for almost every $\theta$,
\[\lim_{n \rightarrow \infty} \frac{a_1+a_1+\ldots+a_n - \max\{a_i : i=1,\ldots,n\}}{n \log n} = \frac{1}{\log 2}.\]  That is to say that the sums grow like $n \log n$ except for rare large partial quotients, corresponding to visits to the cusp in the corresponding flow associated to $\gamma$ in the modular surface $SL_1(\mathbb{R})/SL_2(\mathbb{Z})$ (the paper of Khinchin referenced above contains a proof of convergence \textit{in measure} of the ratio $(a_1+...+a_n)/n \log n$).  This sort of normalization of trimmed sums was extended to a very general setting in \cite{MR1961618}; ergodic sums of nonintegrable functions under exponentially CF-mixing maps, for example our own
\[\rho(\Omega_n)=E(a_1 (\theta_0)) + \cdots + E(a_1 (\theta_{n-1})),\]
may be compared to trimmed sums of regularly varying comparison sequences.

Off the interval $(1/2,1)$, we have the trivial inequality $a_1(\theta)>(1/2)E(a_1 (\theta))>a_1(\theta)/3$.  As the measure preserved by $g$ is mutually absolutely continuous with respect to the Gauss measure, it follows from \cite{MR1961618} that the trimmed sums in our system are almost surely not too far removed from the trimmed sums of partial quotients:
\[\frac{1}{2}\sum_{i=0}^{n-1} E(a_1(\theta_i)) - \frac{1}{2}\max\{E(a_1(\theta_i)):i=0,1,\ldots,n-1\} \sim n \log n.\]
However, we have as a corollary of Theorem \ref{theorem - main result} that
\[\limsup_{n \rightarrow \infty} \frac{\rho_n(x)}{\log n \log \log n} = \infty.\]

\section*{Acknowledgements}
The author is supported by the Center for Advanced Studies at Ben Gurion University and the Israeli Council for Higher Education.

\bibliographystyle{plain}
\bibliography{bibfile}
\end{document}